\documentclass[11pt]{amsart}
\usepackage{mathrsfs}
\usepackage{amssymb,latexsym}
\usepackage{pstricks,color}

 \numberwithin{equation}{section}

\newtheorem{theorem}{Theorem}[section]

\newtheorem{remark}[theorem]{Remark}

\newtheorem{lemma}[theorem]{Lemma}

\newtheorem{example}[theorem]{Example}

\newcommand{\een}{\end{enumerate}}
\newcommand{\blem}{\begin{lem}}
\newcommand{\elem}{\end{lem}}
\newcommand{\bcl}{\begin{clm}}
\newcommand{\ecl}{\end{clm}}
\newcommand{\bthm}{\begin{thm}}
\newcommand{\ethm}{\end{thm}}
\newcommand{\bpr}{\begin{prop}}
\newcommand{\epr}{\end{prop}}
\newcommand{\bco}{\begin{cor}}
\newcommand{\eco}{\end{cor}}
\newcommand{\bcon}{\begin{conj}}
\newcommand{\econ}{\end{conj}}
\newcommand{\bde}{\begin{defn}}
\newcommand{\ede}{\end{defn}}
\newcommand{\bex}{\begin{exa}}
\newcommand{\eexa}{\end{exa}}
\newcommand{\bobs}{\begin{obs}}
\newcommand{\eobs}{\end{obs}}
\newcommand{\bexe}{\begin{exe}}
\newcommand{\eexe}{\end{exe}}

\title{Identities involving Narayana polynomials and Catalan numbers}

\begin{document}
\maketitle
\begin{center}
Toufik Mansour$^\dag$ and Yidong Sun\footnote{Corresponding author:
Yidong Sun, sydmath@yahoo.com.cn.}$^\ddag$

$^\dag$Department of Mathematics, University of Haifa, 31905 Haifa, Israel\\
$^\ddag$Department of Mathematics, Dalian Maritime University, 116026 Dalian, P.R. China\\[5pt]

{\it $^\dag$toufik@math.haifa.ac.il, $^\ddag$sydmath@yahoo.com.cn}
\end{center}\vskip0.5cm

\subsection*{Abstract} We first establish the result that the Narayana
polynomials can be represented as the integrals of the Legendre
polynomials. Then we represent the Catalan numbers in terms of the
Narayana polynomials by three different identities. We give three
different proofs for these identities, namely, two algebraic proofs
and one combinatorial proof. Some applications are also given which
lead to many known and new identities.

{\bf Keywords}: Narayana polynomials, Legendre polynomials, Inverse
relations, Catalan numbers

\noindent {\sc 2000 Mathematics Subject Classification}: Primary
05A05, 05A15

\section{Introduction}
The {\em Catalan numbers} \cite[Sequence A000108]{Slo} are defined
by $C_n=\frac{1}{n+1}\binom{2n}{n}$, for all $n\geq0$.  The {\em
Narayana polynomials} $\mathfrak{N}_n(q)$ and the {\em associated
Narayana polynomials} $\mathcal{N}_n(q)$ \cite{bonin} are defined by
\begin{eqnarray*}
\mathfrak{N}_n(q)=\sum_{k=1}^n\frac{1}{n}\binom{n}{k-1}\binom{n}{k}q^{k}\mbox{
and }\mathcal{N}_n(q)=q^n\mathfrak{N}_n(q^{-1})=\mathfrak{N}_n(q)/q,
\end{eqnarray*}
for $n\geq 1$, with the the initial values
$\mathfrak{N}_0(q)=\mathcal{N}_0(q)=1$. The coefficients
$N_{n,k}=\frac{1}{n}\binom{n}{k-1}\binom{n}{k}$ with $N_{0,0}=1$ are
called the {\em Narayana\ numbers}, and it is well known that the
sequence $\{\mathfrak{N}_n(1)\}_{n\geq 0}$ is the sequence of the
Catalan numbers, while the sequence $\{\mathfrak{N}_n(2)\}_{n\geq
0}$ is the sequence of the large Schr\"{o}der numbers \cite[Sequence
A006318]{Slo}. The Narayana polynomials and associated Narayana
polynomials have been considered by several authors, see
\cite{bonin,coker,rogers,simull,rogershap,sulanke}. For instance,
Bonin, Shapiro and Simion \cite{bonin} showed that the polynomial
$\mathfrak{N}_n(1+q)$ is a $q$-analog of the $n$-th large Schr\"oder
numbers. Coker \cite{coker} provided several different expressions:
\begin{eqnarray}
\sum_{k=1}^n\frac{1}{n}\binom{n}{k-1}\binom{n}{k}q^{k-1}
&=&\sum_{k=0}^{[\frac{n-1}{2}]}\binom{n-1}{2k}C_kq^{k}(1+q)^{n-2k-1}, \label{eqn a1}\\
\sum_{k=1}^n\frac{1}{n}\binom{n}{k-1}\binom{n}{k}q^{2(k-1)}(1+q)^{2(n-k)}
&=&\sum_{k=0}^{n-1}\binom{n-1}{k}C_{k+1}q^{k}(1+q)^{k}.\label{eqn
b1}
\end{eqnarray}
Identity \eqref{eqn a1} was studied by Simion and Ullman
\cite{simull} and proved combinatorially by Chen, Deng and Du
\cite{chendeng}. Later, Chen, Yan and Yang \cite{chenyan} proved
\eqref{eqn a1} and \eqref{eqn b1} combinatorially in terms of
weighted $2$-Motzkin paths. Recently, the authors \cite{mansun1}
presented a new expression for Narayana polynomials
\begin{eqnarray}\label{eqn c1}
\mathfrak{N}_n(q)=\sum_{k=0}^n\frac{1}{n}\binom{n}{k-1}\binom{n}{k}q^k
=\frac{1}{n+1}\sum_{k=0}^{n}\binom{n+1}{k}\binom{2n-k}{n}(q-1)^k,
\end{eqnarray}
whose generalized version is appeared in \cite[Example
2.13]{mansun2}. Chen and Pang \cite{chenpeng} also independently
deduced (\ref{eqn c1}) combinatorially.

The main result of this paper can be formulated as follows.

\begin{theorem}\label{thmm}
For any integer $n\geq 0$, there hold
\begin{eqnarray}
C_n&=&\sum_{k=0}^{n}\frac{2k+1}{2n+1}
\binom{2n+1}{n-k}\mathfrak{N}_k(q)(1-q)^{n-k}, \label{eqn 3.7}\\
q^{\frac{n}{2}+1}C_{\frac{n}{2}}&=&\sum_{k=0}^{n}(-1)^{n-k}\binom{n}{k}\mathfrak{N}_{k+1}(q)(1+q)^{n-k},\label{eqn 3.8}\\
q^{n+2}C_{n+1}&=&\sum_{k=0}^{n}(-1)^{n-k}\binom{n}{k}\mathfrak{N}_{k+1}(q^2)(1-q)^{2(n-k)},\label{eqn
3.9}
\end{eqnarray}
where $C_{\frac{n}{2}}$ is zero if $n$ is odd.
\end{theorem}

In this paper, we first establish the result that the Narayana
polynomials can be represented as the integrals of the Legendre
polynomials \cite{comtet} in Sections~\ref{sec0}. Then we give three
different proofs for Theorem \ref{thmm}, see
Sections~\ref{sec1}-\ref{sec3}, including two algebraical proofs and
one combinatorial proof. Some applications are also given which lead
to many known and new identities.

\section{Narayana polynomials}\label{sec0}
Recall that the Legendre polynomials $P_n(x)$ \cite{comtet,riordan},
which are most familiar in the form
\begin{eqnarray*}
P_n(x)=2^{-n}\sum_{k=0}^{[\frac{n}{2}]}(-1)^k\binom{n-k}{k}\binom{2n-2k}{n-k}x^{n-2k},
\end{eqnarray*}
have an alternate expression, namely,
\begin{eqnarray*}
P_n(x)=\sum_{k=0}^{n}\binom{n}{k}\binom{n+k}{k}\Big(\frac{x-1}{2}\Big)^k=
\sum_{k=0}^{n}\binom{n+k}{n-k}\binom{2k}{k}\Big(\frac{x-1}{2}\Big)^k,
\end{eqnarray*}
so that
\begin{eqnarray}\label{eqn a2}
P_n(2x-1)=\sum_{k=0}^{n}\binom{n+k}{n-k}\binom{2k}{k}(x-1)^k.
\end{eqnarray}

Note that an equivalent form of (\ref{eqn c1}) is
\begin{eqnarray}\label{eqn b2}
\sum_{k=0}^n\frac{1}{n}\binom{n}{k-1}\binom{n}{k}q^k
&=&\sum_{k=0}^{n}\binom{n+k}{n-k}\frac{1}{k+1}\binom{2k}{k}(q-1)^{n-k}.
\end{eqnarray}
Then (\ref{eqn a2}) and (\ref{eqn b2}) generate the following
result.
\begin{theorem}\label{theo 2.1}
For any integer $n\geq 1$, there holds
\begin{eqnarray*}
\mathfrak{N}_n(q)&=&(q-1)^{n+1}\int_{0}^{\frac{q}{q-1}}P_n(2x-1)dx\\
                 &=&q(q-1)^{n}\int_{0}^{1}P_n(\frac{2q}{q-1}x-1)dx.
\end{eqnarray*}
\end{theorem}
\begin{proof} It is clear that $\mathfrak{N}_n(0)=0$ for all $n\geq 1$. Then we
have
\begin{eqnarray*}
\lefteqn{q(q-1)^{n}\int_{0}^{1}P_n(\frac{2q}{q-1}x-1)dx}\\
&=&(q-1)^{n+1}\int_{0}^{\frac{q}{q-1}}P_n(2x-1)dx\\
&=&(q-1)^{n+1}\sum_{k=0}^{n}\binom{n+k}{n-k}\binom{2k}{k}\int_{0}^{\frac{q}{q-1}}(x-1)^kdx\\
&=&\sum_{k=0}^{n}\binom{n+k}{n-k}\frac{1}{k+1}\binom{2k}{k}(q-1)^{n-k}
\\
&+&(q-1)^{n+1}\sum_{k=0}^{n}\binom{n+k}{n-k}\frac{1}{k+1}\binom{2k}{k}(-1)^{k}\\
&=&\mathfrak{N}_n(q)-\mathfrak{N}_n(0)(1-q)^{n+1}=\mathfrak{N}_n(q),
\end{eqnarray*}
which completes the proof.
\end{proof}

Theorem \ref{theo 2.1} signifies that many classical sequences such
as Catalan numbers and Schr\"{o}der numbers can be represented as
the integrals of Legendre polynomials.

\begin{example}\label{example a}
$(i)$ Let $q=-1$. Using the parity identity \cite{coker,sun}
\begin{eqnarray}
\mathfrak{N}_n(-1)=\left\{
\begin{array}{ll}
0 & {\rm if}\ n=2r,\\[5pt]
(-1)^{r+1}C_r & {\rm if}\ n=2r+1,
\end{array}\right.\label{parity}
\end{eqnarray}
we have for $n\geq 0$,
\begin{eqnarray*}
C_{n}=2^{2n+1}\int_{0}^{1}P_{2n+1}(x-1)dx.
\end{eqnarray*}
$(ii)$ Let $q=2$, we have that the large Schr\"{o}der numbers
$\mathfrak{N}_n(2)$ satisfy
\begin{eqnarray*}
\mathfrak{N}_n(2)=\int_{0}^{2}P_{n}(2x-1)dx.
\end{eqnarray*}
\end{example}

\begin{remark} Simons \cite{simons} established the following curious (in
fact, it is not curious) identity
\begin{equation*}
\sum_{k=0}^n \frac{(-1)^{n-k}(n+k)!(1+x)^k}{(n-k)!k!^2}
=\sum_{k=0}^n \frac{(n+k)!x^k}{(n-k)!k!^2},
\end{equation*}
or equivalently
\begin{equation}\label{eqn aa}
\sum_{k=0}^n(-1)^{n-k}\binom{n+k}{n-k}\binom{2k}{k}(1+x)^k
=\sum_{k=0}^n\binom{n+k}{n-k}\binom{2k}{k}x^k,
\end{equation}
which was proved by Chapman \cite{chapman}, Prodinger
\cite{prodinger}, Wang and Sun \cite{wangsun}. It has been pointed
out by Hirschhorn \cite{hirschhorn} that \eqref{eqn aa} is a special
case of the Pfaff identity \cite{gasper}. Recently, Munarini
\cite{munarini} gave a generalization of \eqref{eqn aa}.

Obviously, (\ref{eqn a2}) and (\ref{eqn aa}) generate that
\begin{eqnarray*}\label{eqn ab}
(-1)^{n}P_{n}(-2x-1)=P_{n}(2x+1),
\end{eqnarray*}
which can be easily derived by the generating function of Legendre
polynomials \cite{comtet},
\begin{eqnarray*}
\sum_{n\geq 0}P_n(x)t^n=\frac{1}{\sqrt{1-2xt+t^2}}.
\end{eqnarray*}
\end{remark}



\section{Proof of Theorem~\ref{thmm} and inverse
relations}\label{sec1} In this section, using three well known
inverses relations, we present our first proof for
Theorem~\ref{thmm}. The {\em Legendre inverse relation} reads
\cite{riordan}
\begin{eqnarray}\label{LIR}
A_n=\sum_{k=0}^n\binom{n+k}{n-k}B_k \hskip0.5cm \Longleftrightarrow
\hskip0.5cm
B_n=\sum_{k=0}^{n}(-1)^{n-k}\frac{2k+1}{2n+1}\binom{2n+1}{n-k}A_k,
\end{eqnarray}
and the {\em left-inversion formula} \cite{corsani} reads
\begin{eqnarray}\label{LIF}
A_n=\sum_{k=0}^{[\frac{n}{s}]}\binom{n+p}{sk+p}B_k \hskip0.5cm
\Longrightarrow \hskip0.5cm
B_n=\sum_{k=0}^{sn}(-1)^{sn-k}\binom{sn+p}{k+p}A_k,
\end{eqnarray}
which, in the case $s=1, p=0$, implies the {\em binomial inverse
relation}
\begin{eqnarray}\label{BIR}
A_n=\sum_{k=0}^{n}\binom{n}{k}B_k \hskip0.5cm \Longleftrightarrow
\hskip0.5cm B_n=\sum_{k=0}^{n}(-1)^{n-k}\binom{n}{k}A_k.
\end{eqnarray}

Now we are ready to present the proof of Theorem~\ref{thmm}.

\subsection{Proof of \eqref{eqn 3.7}}
Rewriting (\ref{eqn b2}), we have
\begin{eqnarray*}
\frac{\mathfrak{N}_n(q)}{(q-1)^{n}}
&=&\sum_{k=0}^{n}\binom{n+k}{n-k}C_k(q-1)^{-k},
\end{eqnarray*}
and using \eqref{LIR}, we obtain two expressions for the Catalan
numbers,
\begin{eqnarray*}
C_n&=&\sum_{k=0}^{n}(-1)^{n-k}\frac{2k+1}{2n+1}
\binom{2n+1}{n-k}\mathfrak{N}_k(q)(q-1)^{n-k},
\end{eqnarray*}
which completes the proof of \eqref{eqn 3.7}.

\subsection{Proof of \eqref{eqn 3.8}}
Rewriting (\ref{eqn a1}) in another form after replacing
$n$ by $n+1$,
\begin{eqnarray*}
\frac{\mathfrak{N}_{n+1}(q)}{(1+q)^{n}}
=\sum_{k=0}^{[\frac{n}{2}]}\binom{n}{2k}C_kq^{k+1}(1+q)^{-2k},
\end{eqnarray*}
and using \eqref{LIF} in the case $s=2,p=0$, we deduce another
expression for Catalan numbers,
\begin{eqnarray}\label{eqn catlan2}
q^{n+1}C_n=\sum_{k=0}^{2n}(-1)^{k}\binom{2n}{k}\mathfrak{N}_{k+1}(q)(1+q)^{2n-k},
\end{eqnarray}
which motivates us to consider the following related summation
\begin{eqnarray}\label{eqn catlan3}
f_n(q)=\sum_{i=1}^{2n+2}f_iq^i=\sum_{k=0}^{2n+1}(-1)^{k}\binom{2n+1}{k}\mathfrak{N}_{k+1}(q)(1+q)^{2n+1-k}.
\end{eqnarray}
\begin{lemma}\label{lemma 1}
For all $n\geq0$, $f_n(q)=0$.
\end{lemma}
\begin{proof}
Comparing the coefficients of two sides in (\ref{eqn catlan3}), we
have
\begin{eqnarray*}
f_m=\sum_{j=0}^m\sum_{k=0}^{2n+1}(-1)^{k}\binom{2n+1}{k}N_{k+1,j}\binom{2n+1-k}{m-j}.
\end{eqnarray*}
Noting that $N_{k+1,j}\binom{2n+1-k}{m-j}$ is a polynomial on $k$
with degree $m+j-2$, which does not exceed $2n$ when $1\leq m\leq
n+1$. According to the well-known difference formula
\begin{eqnarray*}
\sum_{k=0}^{n}(-1)^{k}\binom{n}{k}(x-k)^r=\left\{
\begin{array}{ll}
0 & {\rm if}\ 0\leq r<n,\\[5pt]
n! & {\rm if}\ r=n,
\end{array}\right.
\end{eqnarray*}
we can derive that each inner sum is zero in $f_m$ for $1\leq m\leq
n+1$. Note that $f_n(q)=q^{2n+3}f_n(q^{-1})$ by
$q^{n+1}\mathfrak{N}_n(q^{-1})=\mathfrak{N}_n(q)$, which implies
that $f_m=0$ for $n+2\leq m\leq 2n+2$. Hence,  $f_n(q)=0$ for $n\geq
0$, as claimed.
\end{proof}

By combining \eqref{eqn catlan2} and \eqref{eqn catlan3}, using
Lemma \ref{lemma 1}, we obtain \eqref{eqn 3.8}.

\subsection{Proof of \eqref{eqn 3.9}}
Rewriting \eqref{eqn b1} in another form after replacing $n$ by
$n+1$,
\begin{eqnarray*}
\mathfrak{N}_{n+1}\left(\frac{q^2}{(1+q)^2}\right)}{(1+q)^{2n+2}
=\sum_{k=0}^{n}\binom{n}{k}C_{k+1}q^{k+2}(1+q)^{k},
\end{eqnarray*}
and using \eqref{BIR}, we deduce that
\begin{eqnarray*}
q^{n+2}(1+q)^nC_{n+1}=\sum_{k=0}^{n}(-1)^{n-k}\binom{n}{k}\mathfrak{N}_{k+1}\left(\frac{q^2}{(1+q)^2}\right)(1+q)^{2k+2}.
\end{eqnarray*}
Replacing $q$ by $\frac{q}{1-q}$, after simplification, we get
\eqref{eqn 3.9}.

\subsection{Applications}
Theorem \ref{thmm} can produce numerical known or new identities.
For instance,

\begin{itemize}
\item The case $q=-1$ in \eqref{eqn 3.7} together with \eqref{parity}, lead to a new identity
\begin{eqnarray*}
(2^n-1)C_n=\sum_{r=0}^{[\frac{n-1}{2}]}(-1)^r\frac{4r+3}{2n+1}\binom{2n+1}{n-2r-1}2^{n-2r-1}C_r,
\end{eqnarray*}
which implies that $C_{2k}\equiv0\ mod\ 2$ and $C_{2k-1}\equiv
C_{k-1}\ mod\ 2$ for $k\geq 1$, from which one can easily derive
that $C_{n}$ is odd if and only if $n=2^k-1$ for some $k\geq 0$.

\item Taking the coefficient of $q^n$ in both sides of \eqref{eqn
3.7}, we get another parity identity
\begin{eqnarray}\label{eqn 3.10}
\sum_{k=0}^{n}(-1)^{k}\frac{2k+1}{2n+1} \binom{2n+1}{n-k}=0,
\hskip0.6cm (n\geq 1),
\end{eqnarray}
which has been proved by Chen, Li and Shapiro \cite{chenli}.

\item The case $q=1$ in \eqref{eqn 3.8} leads to a new identity
\begin{eqnarray*}
C_{n}=\sum_{k=0}^{2n}(-1)^{k}\binom{2n}{k}C_{k+1}2^{2n-k}.
\end{eqnarray*}

\item The case $q=-1$ in \eqref{eqn 3.9} leads to a known identity
\cite{aigner,coker,riordan}
\begin{eqnarray*}
C_{n+1}=\sum_{k=0}^{n}(-1)^{k}\binom{n}{k}C_{k+1}4^{n-k}.
\end{eqnarray*}
and the case $q=\sqrt{-1}$ in \eqref{eqn 3.9} leads to the Touchard
identity \cite{aigner,coker}
\begin{eqnarray*}
C_{n+1}=\sum_{k=0}^{n}\binom{n}{2k}C_{k}2^{n-2k}.
\end{eqnarray*}

\item Let $q=\sqrt{2}$ in (\ref{eqn 3.9}), by the relation
$(1-\sqrt{2})^n=(P_n+P_{n-1})-P_n\sqrt{2}$, where $P_n$ is the
$n$-th Pell number (defined by the recurrence relation
$P_{n+1}=2P_{n}+P_{n-1}$ with $P_{-1}=1,P_0=0$), we have new
identities involving Catalan numbers, large Schr\"{o}der numbers,
and Pell numbers
\begin{eqnarray*}
2^{n+1}C_{2n+1}&=&\sum_{k=0}^{2n}(-1)^{k}\binom{2n}{k}\mathfrak{N}_{k+1}(2)P_{4n-2k-1},\\
2^{n+1}C_{2n+2}&=&\sum_{k=0}^{2n+1}(-1)^{k}\binom{2n+1}{k}\mathfrak{N}_{k+1}(2)P_{4n-2k+2}.
\end{eqnarray*}

\item Let $q=\sqrt{5}$ in (\ref{eqn 3.9}), by the relation
$(\frac{1-\sqrt{5}}{2})^{n+1}=\frac{L_n-F_n\sqrt{5}}{2}$, where
$L_n$ and $F_n$ are respectively the $n$-th Lucas number and the
$n$-th Fibonacci number (defined by the same recurrence relation
$G_{n+1}=G_{n}+G_{n-1}$ with $G_{-1}=2,G_0=1$ for $L_n$ and
$G_{-1}=0,G_0=1$ for $F_n$), we have new identities involving
Catalan numbers, Lucas numbers, and Fibonacci numbers
\begin{eqnarray*}
5^{n+1}C_{2n+1}&=&\sum_{k=0}^{2n}(-1)^{k}\binom{2n}{k}\mathfrak{N}_{k+1}(5)L_{4n-2k-1}2^{4n-2k-1},\\
5^{n+1}C_{2n+2}&=&\sum_{k=0}^{2n+1}(-1)^{k}\binom{2n+1}{k}\mathfrak{N}_{k+1}(5)F_{4n-2k+1}2^{4n-2k+1}.
\end{eqnarray*}
\end{itemize}

\section{Proof of Theorem \ref{thmm} and generating
functions}\label{sec2} In this section we present our second proof
for Theorem~\ref{thmm} which is based on generating function
techniques.

Recall that $C(x)=\frac{1-\sqrt{1-4x}}{2x}$ is the generating
function for the Catalan numbers $C_n=\frac{1}{n+1}\binom{2n}{n}$,
which satisfies the relation $C(x)=1+xC(x)^2=\frac{1}{1-xC(x)}$. By
Lagrange inversion formula~\cite{wilf}, one can deduce that
\begin{eqnarray}\label{eqn 4.1}
[x^{n-k}]C(x)^{2k+1}=\frac{2k+1}{2n+1} \binom{2n+1}{n-k}.
\end{eqnarray}
Define $\Omega(q,x)=\sum_{n\geq 0}\mathfrak{N}_{n}(q)x^n$, then
$\Omega(q,x)$ has the explicit expression \cite{deutsch}
\begin{eqnarray*}
\Omega(q,x)=\frac{1+x-qx-\sqrt{1-2x+x^2-2qx-2qx^2+q^2x^2}}{2x},
\end{eqnarray*}
which can be rewritten as
\begin{eqnarray*}
\Omega(q,x)=\frac{1}{1+x-qx}C\left(\frac{x}{(1+x-qx)^2}\right)=1+\frac{qx}{1-x-qx}C\left(\frac{qx^2}{(1-x-qx)^2}\right).
\end{eqnarray*}
By the Cauchy Residue Theorem, we have
\begin{eqnarray*} \lefteqn{\sum_{k=0}^{n}\frac{2k+1}{2n+1}
\binom{2n+1}{n-k}\mathfrak{N}_k(q)(1-q)^{n-k}}\\
&=&\sum_{k=0}^{n}\underset{x}{Res}\frac{C((1-q)x)^{2k+1}}{x^{n-k+1}}\underset{y}{Res}\frac{\Omega(q,y)}{y^{k+1}}\\
&=&\underset{x}{Res}\frac{C((1-q)x)}{x^{n+1}}\sum_{k=0}^{n}x^kC((1-q)x)^{2k}\underset{y}{Res}\frac{\Omega(q,y)}{y^{k+1}}\\
&=&\underset{x}{Res}\frac{C((1-q)x)}{x^{n+1}}\Omega(q,xC((1-q)x)^2)\\
&=&\underset{x}{Res}\frac{C((1-q)x)}{1+(1-q)xC((1-q)x)^2}C\left(\frac{xC((1-q)x)^2}{(1+(1-q)xC((1-q)x)^2)^2}\right){x^{-n-1}}\\
&=&\underset{x}{Res}{C(x)}{x^{-n-1}}=C_n,
\end{eqnarray*}
\begin{eqnarray*}
\lefteqn{\sum_{k=0}^{n}(-1)^{n-k}\binom{n}{k}\mathfrak{N}_{k+1}(q)(1+q)^{n-k}}\\
&=&\underset{x}{Res}\frac{(1-(1+q)x)^n(\Omega(q,x)-1)}{x^{n+2}}    \\
&=&\underset{z}{Res}\frac{\Omega(q,\frac{z}{1+(1+q)z})-1}{z^{n+2}},\hskip2cm ({\rm by\ setting}\  x=\frac{z}{1+(1+q)z})\\
&=&\underset{z}{Res}\left\{\frac{qx}{1-(1+q)x}C\left(\frac{qx^2}{(1-(1+q)x)^2}\right)\right\}_{x=\frac{z}{1+(1+q)z}}z^{-n-2}\\
&=&\underset{z}{Res}\frac{qC(qz^2)}{z^{n+1}}=q^{\frac{n}{2}+1}C_{\frac{n}{2}},
\end{eqnarray*}
and
\begin{eqnarray*}
\lefteqn{\sum_{k=0}^{n}(-1)^{n-k}\binom{n}{k}\mathfrak{N}_{k+1}(q^2)(1-q)^{2(n-k)}}\\
&=&\underset{x}{Res}\frac{(1-(1-q)^2x)^n(\Omega(q^2,x)-1)}{x^{n+2}}    \\
&=&\underset{z}{Res}\frac{\Omega(q^2,\frac{z}{1+(1-q)^2z})-1}{z^{n+2}},\hskip2cm ({\rm by\ setting}\  x=\frac{z}{1+(1-q)^2z})\\
&=&\underset{z}{Res}\left\{\frac{q^2x}{1-(1+q^2)x}C\left(\frac{q^2x^2}{(1-(1+q^2)x)^2}\right)\right\}_{x=\frac{z}{1+(1-q^2)z}}z^{-n-2}     \\
&=&\underset{z}{Res}\frac{q^2}{1-2qz}C\left(\frac{q^2z^2}{(1-2qz)^2}\right){z^{-n-1}}\\
&=&\underset{w}{Res}\frac{q^{n+2}(1+2w)^nC(w^2)}{w^{n+1}},\hskip1.8cm ({\rm by\ setting}\  z=\frac{w}{q(1+2w)})   \\
&=&q^{n+2}\sum_{k=0}^{n}\binom{n}{2k}C_k2^{n-2k}=q^{n+2}C_{n+1},
\end{eqnarray*}
which give \eqref{eqn 3.7}, \eqref{eqn 3.8} and \eqref{eqn 3.9},
respectively. Note that the last equation follows by the well-known
Touchard's identity, which can also be derived by setting $q=1$ in
(\ref{eqn a1}) after replacing $n$ by $n+1$.

\section{Combinatorial Proof of Theorem \ref{thmm}}\label{sec3}

In order to give the combinatorial proof of (\ref{eqn 3.7}), we need
the following definitions. A {\em Dyck path} of length $2n$ is a
lattice path from $(0,0)$ to $(2n,0)$ in the first quadrant of
$xoy$-plane, consisting of up-steps $u=(1,1)$ and down-steps
$d=(1,-1)$, which never passes below the $x$-axis. We will refer to
$n$ as the semilength of the path. It is well known that the set of
Dyck paths of semilength $n$ is counted by the Catalan number
$C_n=\frac{1}{2n+1}\binom{2n+1}{n}$. A {\em peak} in a Dyck path is
an occurrence of $ud$. By the {\em hight} of a step we mean the
ordinate of its endpoint. By a {\em return step} we mean a down-step
ending at height zero. Dyck paths that have exactly one return step
are said to be {\em primitive}. If $D_1$ and $D_2$ are Dyck paths,
we define $D_1D_2$ to be the {\em concatenation} of $D_1$ and $D_2$.

A {\em weighted Dyck path} is a Dyck path $D$ for which every
up-step is endowed with a weight. The weight of a Dyck path $D$ is
the product of the weights of its up-steps, the weight of a set $S$
of Dyck paths means the sum of the weights of $D$ in $S$.

Let $\mathscr{D}_{n,k}$ denote the set of weighted Dyck paths of
length $2n$ consisting of an ordered sequence of $2k+2$ Dyck paths
$D^{(0)}, D^{(1)}, \dots, D^{(2k+1)}$ such that
\begin{itemize}
\item $D^{(0)}$ is a Dyck path of length $2k$ for $0\leq k\leq n$
with an up-step in each peak weighted by $q$ and other up-steps
weighted by $1$;
\item There are totally $n-k$ up-steps in the rest $2k+1$ Dyck
paths, and all up-steps of each $D^{(i)}$ are weighted by 1 or $-q$
for $1\leq i\leq 2k+1$, i.e., such up-steps can be regarded to be
weighted by $(1-q)$;
\item Each $D^{(i)}$ is inserted into the $i$-th endpoint of
$D^{(0)}$ including the beginning point for $1\leq i\leq 2k+1$.
\end{itemize}
Let $\mathscr{D}_n=\bigcup_{k=0}^n\mathscr{D}_{n,k}$. For any $D\in
\mathscr{D}_n$, denote by $w(D)$ the weight of $D$. Let
$\tilde{\mathscr{D}}_n$ denote the subset of Dyck paths
$D\in\mathscr{D}_n$ such that all up-steps in $D$ are weighted by 1,
such Dyck paths only appear in $\mathscr{D}_{n,0}$.

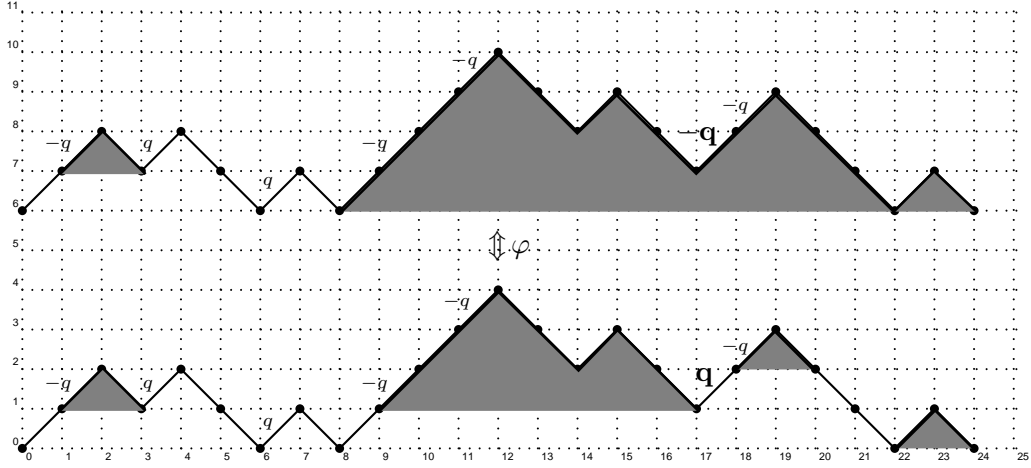
\begin{figure}[h] \setlength{\unitlength}{0.5mm}
\begin{center}
\begin{pspicture}(13,5.9)
\psset{xunit=15pt,yunit=15pt}\psgrid[subgriddiv=1,griddots=4,
gridlabels=4pt](0,0)(25,11)

\psline(0,0)(2,2)(3,1)(4,2)(6,0)(7,1)(8,0)(12,4)(14,2)(15,3)(17,1)(19,3)(22,0)(23,1)(24,0)

\pscircle*(0,0){0.06}\pscircle*(1,1){0.06}\pscircle*(2,2){0.06}
\pscircle*(3,1){0.06}\pscircle*(4,2){0.06}\pscircle*(5,1){0.06}
\pscircle*(6,0){0.06}\pscircle*(7,1){0.06}\pscircle*(8,0){0.06}
\pscircle*(9,1){0.06}\pscircle*(10,2){0.06}\pscircle*(11,3){0.06}
\pscircle*(12,4){0.06}\pscircle*(13,3){0.06}\pscircle*(14,2){0.06}
\pscircle*(15,3){0.06}\pscircle*(16,2){0.06}\pscircle*(17,1){0.06}
\pscircle*(18,2){0.06}\pscircle*(19,3){0.06}\pscircle*(20,2){0.06}
\pscircle*(21,1){0.06}\pscircle*(22,0){0.06}\pscircle*(23,1){0.06}
\pscircle*(24,0){0.06}

\put(3.2,.3){\tiny$q$}

\put(.3,.8){\tiny$-q$}\put(1.6,.8){\tiny$q$}
\put(4.5,.8){\tiny$-q$}\put(8.95,.9){$\bf q$}

\put(5.6,1.9){\tiny$-q$}

\put(9.3,1.3){\tiny$-q$}

{\put(.5,.5){\psline[unit=18pt,linewidth=1pt,fillstyle=solid,fillcolor=gray](0,0)(.9,.9)(1.8,0)}}
{\put(11.63,0){\psline[unit=18pt,linewidth=1pt,fillstyle=solid,fillcolor=gray](0,0)(.8,.8)(1.6,0)}}
{\put(9.5,1.05){\psline[unit=18pt,linewidth=1pt,fillstyle=solid,fillcolor=gray](0,0)(.8,.8)(1.6,0)}}
{\put(4.75,.5){\psline[unit=18pt,linewidth=1pt,fillstyle=solid,fillcolor=gray](0,0)(2.5,2.5)(4.18,.85)(5.,1.7)(6.7,0)}}

\put(6.2,2.6){$\Updownarrow$}\put(6.5,2.6){$\varphi$}

\psline(0,6)(2,8)(3,7)(4,8)(6,6)(7,7)(8,6)(12,10)(14,8)(15,9)(17,7)(19,9)(22,6)(23,7)(24,6)

\pscircle*(0,6){0.06}\pscircle*(1,7){0.06}\pscircle*(2,8){0.06}
\pscircle*(3,7){0.06}\pscircle*(4,8){0.06}\pscircle*(5,7){0.06}
\pscircle*(6,6){0.06}\pscircle*(7,7){0.06}\pscircle*(8,6){0.06}
\pscircle*(9,7){0.06}\pscircle*(10,8){0.06}\pscircle*(11,9){0.06}
\pscircle*(12,10){0.06}\pscircle*(13,9){0.06}\pscircle*(14,8){0.06}
\pscircle*(15,9){0.06}\pscircle*(16,8){0.06}\pscircle*(17,7){0.06}
\pscircle*(18,8){0.06}\pscircle*(19,9){0.06}\pscircle*(20,8){0.06}
\pscircle*(21,7){0.06}\pscircle*(22,6){0.06}\pscircle*(23,7){0.06}
\pscircle*(24,6){0.06}

\put(3.2,3.5){\tiny$q$}

\put(.3,4){\tiny$-q$}\put(1.6,4){\tiny$q$}
\put(4.5,4){\tiny$-q$}\put(8.7,4.1){$\bf -q$}

\put(5.7,5.1){\tiny$-q$}

\put(9.3,4.5){\tiny$-q$}

{\put(.5,3.65){\psline[unit=18pt,linewidth=1pt,fillstyle=solid,fillcolor=gray](0,0)(.9,.9)(1.8,0)}}
{\put(4.25,3.15){\psline[unit=18pt,linewidth=1pt,fillstyle=solid,fillcolor=gray]
(0,0)(3.3,3.3)(4.95,1.65)(5.78,2.45)(7.45,.8)(9.12,2.45)(11.6,0)(12.45,0.85)(13.34,0)}}

\end{pspicture}
\caption{ The involution $\varphi$ on $\mathscr{D}_n\setminus
\tilde{\mathscr{D}}_n$, where the weight $1$ of up-steps is
unlabeled. }\label{fDD0}
\end{center}
\end{figure}

\begin{theorem}\label{theo 4.0}
There exists a sign reversing involution $\varphi$ on the set
$\mathscr{D}_n\setminus \tilde{\mathscr{D}}_n$.
\end{theorem}
\begin{proof}
For any $D\in\mathscr{D}_n\setminus \tilde{\mathscr{D}}_n$, it can
be uniquely written as $D=D_1D_2\cdots D_m$ for some $1\leq m\leq
n$, where $D_i$'s are weighted primitive Dyck paths. Obviously,
there exist at least a $D_i$ such that $D_i$ has an up-step weighted
by $q$ or $-q$. Now we can recursively construct the involution
$\varphi$ as follows. First find the maximum $i$ for $1\leq i\leq m$
such that $D_i$ has an up-step weighted by $q$ or $-q$, then define
$\varphi(D)=D_1\cdots D_{i-1}\varphi(D_i)D_{i+1}\cdots D_m$. Note
that $D_i=uD_i^{*}d$, where $D_i^{*}\in \mathscr{D}_j$ for some
$0\leq j\leq n-1$.
\begin{itemize}
\item  If the first up-step in $D_i=uD_i^{*}d$ has
weight $q$ or $-q$, then $D_i^{*}$ has no up-steps with weight $q$,
otherwise it will contradict with the definition of $D\in
\mathscr{D}_n$. Let $D_i^{'}$ denote the weighted Dyck path obtained
from $D_i$ by changing the sign of the weight of the first up-step.
Then define $\varphi(D_i)=D_i^{'}$;
\item If the first up-step in $D_i=uD_i^{*}d$ has weight 1, then
$D_i^{*}\in \mathscr{D}_j\setminus \tilde{\mathscr{D}}_j$ for some
$1\leq j\leq n-1$, define $\varphi(D_i)=u\varphi(D_i^{*})d$.
\end{itemize}
It is clear that the $\varphi$ is a sign reversing involution on the
set $\mathscr{D}_n\setminus \tilde{\mathscr{D}}_n$. See Figure
\ref{fDD0} for an illustration.
\end{proof}

\begin{proof}[Proof of \eqref{eqn
3.7}] It is clear that the weight of $\tilde{\mathscr{D}}_n$ is the
$n$-th Catalan number $C_n$. According to the definition of
$D\in\mathscr{D}_{n,k}$, it is easy to derive the weight of
$\mathscr{D}_{n,k}$. On one hand, it is well known that the total
weight for $D^{(0)}$ is the Narayana polynomial $\mathfrak{N}_k(q)$.
On the other hand, the total product of the weights of $D^{(1)},
D^{(2)}, \dots, D^{(2k+1)}$ is just the coefficient of $x^{n-k}$ in
$(1-q)^{n-k}C(x)^{2k+1}$. Then by (\ref{eqn 4.1}) we have
$$w(\mathscr{D}_{n,k})=\mathfrak{N}_k(q)(1-q)^{n-k}[x^{n-k}]C(x)^{2k+1}=\frac{2k+1}{2n+1}
\binom{2n+1}{n-k}\mathfrak{N}_k(q)(1-q)^{n-k}.$$

By Theorem \ref{theo 4.0}, we have
$w(\tilde{\mathscr{D}}_n)=\sum_{k=0}^nw(\mathscr{D}_{n,k})$, which
completes the proof.
\end{proof}
\begin{remark}
Specially, let $\bar{\mathscr{D}}_n$ denote the subset of Dyck paths
$D\in\mathscr{D}_n$ such that all up-steps in $D$ are weighted by
$q$ or $-q$, such Dyck paths can only appear in $\mathscr{D}_{n,k}$
for $0\leq k\leq n$ satisfying that $(a)$ $D^{(0)}=(ud)^k$ with
up-steps weighted by $q$, and $(b)$ all up-steps in $D^{(i)}$ are
weighted by $-q$ for $1\leq i\leq k$. Let
$\varphi_{\bar{\mathscr{D}}_n}$ be the $\varphi$ restricted to
$\bar{\mathscr{D}}_n$, it is clear that
$\varphi_{\bar{\mathscr{D}}_n}$ is a sign reversing involution on
$\bar{\mathscr{D}}_n$. Then (\ref{eqn 3.10}) is followed immediately
by $\varphi_{\bar{\mathscr{D}}_n}$ and (\ref{eqn 4.1}).

\end{remark}

In order to give the combinatorial proof of (\ref{eqn 3.8}) and
(\ref{eqn 3.9}), we need the following definitions. A {\em plane
tree } $T$ can be defined recursively (see for example
\cite{stanley}) as a finite set of vertices such that a
distinguished vertex $u$ is called the {\em root} of $T$, and the
remaining vertices are put into an ordered partition
$(T_1,T_2,\cdots,T_m)$ of $m\geq 0$ disjoint non-empty sets, each of
which is a plane tree called the subtree of $u$. The root $u_i$ of
$T_i$ is called the {\em child} of $u$, and $u$ is called the {\em
father} of $u_i$. The {\em out-degree} of a vertex of $T$ is the
number of its subtrees. An {\em internal vertex} of $T$ is a vertex
of out-degree at least one. A vertex of out-degree zero is called a
{\em leaf} of $T$. A {\em complete binary tree} is a plane tree such
that each internal vertex has out-degree two.

A {\em weighted plane tree } is a plane tree for which every vertex
is endowed with a weight. The weight of a plane tree $T$ is the
product of the weights of its vertices, the weight of a set $S$ of
plane trees means the sum of the weights of $T$ in $S$.

Let $\mathscr{P}_{n,k}$ denote the set of weighted plane trees of
$n+2$ vertices such that
\begin{itemize}
\item The leaves have weight $q$;
\item There exist $n-k$ vertices of out-degree one, except for the root, with weight
$-1$ or $-q$, in other words, such vertices can be regarded to be
weighted by $-(1+q)$;
\item All other internal vertices have weight $1$, there may exist
vertices of out-degree one with weight 1.
\end{itemize}
Let $\mathscr{P}_n=\bigcup_{k=0}^n\mathscr{P}_{n,k}$. For any $T\in
\mathscr{P}_n$, denote by $w(T)$ the weight of $T$. Let
$\mathscr{P}_n^*$ denote the subset of $\mathscr{P}_n$ such that
there are at least one vertex, except for the root, of out-degree
one weighted by $1$ or $-1$, let $\tilde{\mathscr{P}}_n$ denote the
subset of $\mathscr{P}_n$ such that the root has out-degree one, and
all other internal vertices have out-degree two. It is clear that
$\mathscr{P}_n^*\cap \tilde{\mathscr{P}}_n=\emptyset$.
\begin{theorem}\label{theo 4.1}
There exists a sign reversing involution $\psi$ on the set
$\mathscr{P}_n\setminus \tilde{\mathscr{P}}_n$.
\end{theorem}
\begin{proof}
Note that a tree $T\in \mathscr{P}_n$ is in $\mathscr{P}_n^*$ if and
only if it contains a non-rooted vertex of out-degree one with
weight either $1$ or $-1$. Consider the first occurrence of such
vertex, denoted by $v$, when traversing the weighted plane tree $T$
in pre-order, (i.e., visiting the root first, then traversing its
subtrees from left to right). Then replace the weight of $v$ in $T$
by $-w(v)$, we obtain another weighted plane tree $T^*$ in
$\mathscr{P}_n^*$, and then define $\psi(T)=T^*$. See Figure
\ref{fDD1} for example.

\begin{figure}[h]
\setlength{\unitlength}{0.5mm}
\begin{center}
\begin{pspicture}(11,4.1)
\psset{xunit=20pt,yunit=20pt}\psgrid[subgriddiv=1,griddots=5,
gridlabels=4pt](0,0)(16,5.5)

\psline(4,5)(2,3.5) \psline(4,5)(4,3.5)(4,2)(4,.5)
\psline(4,5)(6,3.5)
\psline(2,3.5)(1,2)\psline(2,3.5)(3,2)\psline(4,3.5)(6,2)(6,.5)

\pscircle*(4,5){0.06}
\pscircle*(2,3.5){0.06}\pscircle*(4,3.5){0.06}\pscircle*(6,3.5){0.06}
\pscircle*(1,2){0.06}\pscircle*(3,2){0.06}\pscircle*(4,2){0.1}\pscircle*(6,2){0.06}
\pscircle*(4,.5){0.06}\pscircle*(6,.5){0.06}

\put(3.,1.3){$-1$}\put(4.4,1.3){$-1$}
\put(5.3,2){$\rightleftharpoons$}

\psline(12,5)(10,3.5) \psline(12,5)(12,3.5)(12,2)(12,.5)
\psline(12,5)(14,3.5)
\psline(10,3.5)(9,2)\psline(10,3.5)(11,2)\psline(12,3.5)(14,2)(14,.5)

\pscircle*(12,5){0.06}
\pscircle*(10,3.5){0.06}\pscircle*(12,3.5){0.06}\pscircle*(14,3.5){0.06}
\pscircle*(9,2){0.06}\pscircle*(11,2){0.06}\pscircle*(12,2){0.1}\pscircle*(14,2){0.06}
\pscircle*(12,.5){0.06}\pscircle*(14,.5){0.06}

\put(8.6,1.3){$1$}\put(10,1.3){$-1$}
\end{pspicture}
\caption{ The involution on $\mathscr{P}_n^{*}$}\label{fDD1}
\end{center}
\end{figure}
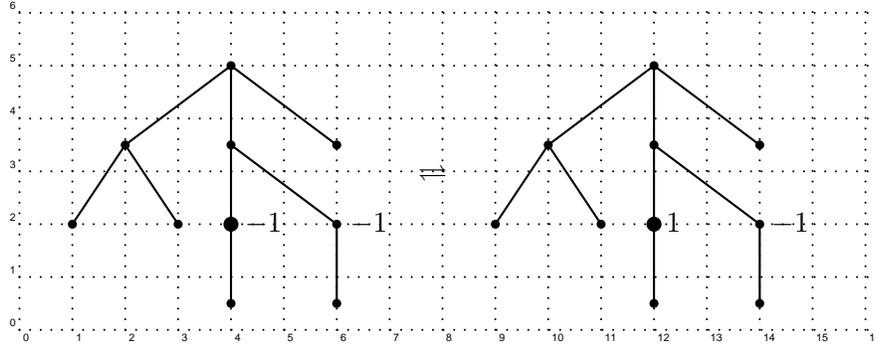

For any $T\in \mathscr{P}_n\setminus (\mathscr{P}_n^*\cup
\tilde{\mathscr{P}}_n)$, we can construct recursively the involution
$\psi$ as follows. First we should consider the following two cases:
\begin{itemize}
\item The root of $T$ has out-degree one, in this case, the unique subtree of the root has
either a vertex of out-degree greater than two or a vertex of
out-degree one with weight $-q$.
\item The root of $T$ has out-degree not less than two.
\end{itemize}
If the root $u$ of $T$ has $i(\geq 2)$ number of subtrees, denoted
by $T_1,T_2,\cdots,T_i$, when $T_1$ is a complete binary tree, then
delete the subtree $T_2$ of the root $u$, regard $T_2$ as the
subtree of the right-most leaf $v$ of $T_1$, and replace the weight
$q$ of $v$ by $-q$, hence we obtain a new weighted plane tree $T^*$
in $\mathscr{P}_n\setminus (\mathscr{P}_n^*\cup
\tilde{\mathscr{P}}_n)$, and define $\psi(T)=T^*$; when $T_1$ is not
a complete binary tree, consider the left-most {\em component} of
$T$, that is the tree $T_1'$ with the root $u$ and $T_1$ as its
unique subtree, then $\psi(T)$ is obtained by adding the subtrees
$T_2,\cdots,T_i$ right to the root $u$ of $\psi(T_1')$ step by step.
See Figure \ref{fDD2} for example.
\begin{figure}[h]
\setlength{\unitlength}{0.5mm}
\begin{center}
\begin{pspicture}(11,4.1)
\psset{xunit=20pt,yunit=20pt}\psgrid[subgriddiv=1,griddots=5,
gridlabels=4pt](0,0)(16,5.5)

\psline(4,5)(2,3.5) \psline(4,5)(4,3.5) \psline(4,5)(6,3.5)(5,2)
\psline(2,3.5)(1,2)\psline(2,3.5)(3,2)\psline(6,3.5)(6,2)(6,.5)

\pscircle*(4,5){0.06}
\pscircle*(2,3.5){0.06}\pscircle*(4,3.5){0.06}\pscircle*(6,3.5){0.06}
\pscircle*(1,2){0.06}\pscircle*(3,2){0.1}\pscircle*(5,2){0.06}\pscircle*(6,2){0.06}
\pscircle*(6,.5){0.06}

\put(2.4,1.3){$q$}\put(4.4,1.3){$-q$}\put(3,2.3){$q$}
\put(5.4,2){$\rightleftharpoons$}

\psline(12,5)(10,3.5) \psline(12,5)(14,3.5)(13,2)
\psline(10,3.5)(9,2)\psline(10,3.5)(11,2)(11,.5)\psline(14,3.5)(14,2)(14,.5)

\pscircle*(12,5){0.06}
\pscircle*(10,3.5){0.06}\pscircle*(11,.5){0.06}\pscircle*(14,3.5){0.06}
\pscircle*(9,2){0.06}\pscircle*(11,2){0.1}\pscircle*(13,2){0.06}\pscircle*(14,2){0.06}
\pscircle*(11,.5){0.06}\pscircle*(14,.5){0.06}

\put(10,1.3){$-q$}\put(7.9,1.3){$-q$}\put(7.9,.3){$q$}
\end{pspicture}
\caption{ The involution on $\mathscr{P}_n\setminus
\tilde{\mathscr{P}}_n$}\label{fDD2}
\end{center}
\end{figure}
If the root $u$ of $T$ has a unique subtree, denoted by $T'$, let
$u'$ be the root of $T'$ which is the only child of $u$.

$\rm (i)$ If the out-degree of $u'$ is greater than two, then
$\psi(T)$ is defined to be the tree $T^*$ in $\mathscr{P}_n\setminus
(\mathscr{P}_n^*\cup \tilde{\mathscr{P}}_n)$ which has the root $u$
with a unique subtree $\psi(T')$; See Figure \ref{fDD3} for example.

\begin{figure}[h]
\setlength{\unitlength}{0.5mm}
\begin{center}
\begin{pspicture}(11,4.8)
\psset{xunit=20pt,yunit=20pt}\psgrid[subgriddiv=1,griddots=5,
gridlabels=4pt](0,0)(16,7)

\psline(4,6.5)(4,5)(2,3.5) \psline(4,5)(4,3.5)
\psline(4,5)(6,3.5)(5,2)
\psline(2,3.5)(1,2)\psline(2,3.5)(2,2)\psline(2,3.5)(3,2)\psline(6,3.5)(6,2)(6,.5)

\pscircle*(4,6.5){0.06}\pscircle*(4,5){0.06}
\pscircle*(2,3.5){0.06}\pscircle*(4,3.5){0.06}\pscircle*(6,3.5){0.06}
\pscircle*(1,2){0.1}\pscircle*(2,2){0.06}\pscircle*(3,2){0.06}\pscircle*(5,2){0.06}\pscircle*(6,2){0.06}
\pscircle*(6,.5){0.06}

\put(.3,1.3){$q$}\put(4.4,1.3){$-q$}\put(1.6,1.3){$q$}
\put(5.4,2.5){$\rightleftharpoons$}

\psline(12,6.5)(12,5)(10,3.5)
\psline(12,5)(14,3.5)(13,2)\psline(12,5)(12,3.5)
\psline(10,3.5)(9,2)(9,.5)\psline(10,3.5)(11,2)\psline(14,3.5)(14,2)(14,.5)

\pscircle*(12,6.5){0.06}\pscircle*(12,5){0.06}
\pscircle*(10,3.5){0.06}\pscircle*(12,3.5){0.06}\pscircle*(14,3.5){0.06}
\pscircle*(9,2){0.1}\pscircle*(11,2){0.06}\pscircle*(13,2){0.06}\pscircle*(14,2){0.06}
\pscircle*(9,.5){0.06}\pscircle*(14,.5){0.06}

\put(10,1.3){$-q$}\put(6.5,1.3){$-q$}\put(6.5,.3){$q$}
\end{pspicture}
\caption{ The involution on $\mathscr{P}_n\setminus
\tilde{\mathscr{P}}_n$}\label{fDD3}
\end{center}
\end{figure}

$\rm (ii)$ If the out-degree of $u'$ is one or two, find the
right-most leaf $v'$,

\begin{itemize}
\item if there exist vertices of weight $-q$ in the path $u'v'$,
then choose the vertex $v$ which is first occurring in the path
$u'v'$, denoted by $T''$ as the subtree of $v$,  if deleting $T''$
in $T'$, the resulting tree is a complete binary tree. Then deleting
the subtree $T''$ in $T$, annexing it to the right of $u$, and
changing the weight $-q$ of $v$ to be $q$, we obtain a new weighted
plane tree $T^*$ in $\mathscr{P}_n\setminus (\mathscr{P}_n^*\cup
\tilde{\mathscr{P}}_n)$, and define $\psi(T)=T^*$; See Figure
\ref{fDD4} for example.
\begin{figure}[h]
\setlength{\unitlength}{0.5mm}
\begin{center}
\begin{pspicture}(11,4.8)
\psset{xunit=20pt,yunit=20pt}\psgrid[subgriddiv=1,griddots=5,
gridlabels=4pt](0,0)(16,7)

\psline(4,6.5)(4,5)(2,3.5) \psline(4,5)(6,3.5)(5,2)
\psline(2,3.5)(1,2)\psline(2,3.5)(3,2)\psline(6,3.5)(6,2)(6,.5)

\pscircle*(4,6.5){0.06}\pscircle*(4,5){0.06}
\pscircle*(2,3.5){0.06}\pscircle*(6,3.5){0.06}
\pscircle*(1,2){0.06}\pscircle*(3,2){0.06}\pscircle*(5,2){0.06}\pscircle*(6,2){0.1}
\pscircle*(6,.5){0.06}

\put(4.4,1.3){$-q$}\put(4.4, .3){$q$}
\put(5.4,2.5){$\rightleftharpoons$}

\psline(12,6.5)(12,5)(10,3.5)\psline(12,6.5)(14,5)
\psline(12,5)(14,3.5)(13,2)
\psline(10,3.5)(9,2)\psline(10,3.5)(11,2)\psline(14,3.5)(14,2)

\pscircle*(12,6.5){0.06}\pscircle*(12,5){0.06}
\pscircle*(10,3.5){0.06}\pscircle*(14,3.5){0.06}
\pscircle*(9,2){0.06}\pscircle*(11,2){0.06}\pscircle*(13,2){0.06}\pscircle*(14,2){0.06}
\pscircle*(14,5){0.1}

\put(10,1.3){$q$}\put(10,3.5){$q$}
\end{pspicture}
\caption{ The involution on $\mathscr{P}_n\setminus
\tilde{\mathscr{P}}_n$}\label{fDD4}
\end{center}
\end{figure}

\item if deleting $T''$ in $T'$, the resulting tree is not a complete
binary tree or if there is no vertex of weight $-q$ in the path
$u'v'$, then $u'$ must have out-degree two. Let $T_1,T_2$ be the
left and right subtrees of $u'$ and $T_1',T_2'$ be the left and
right components of $u'$ respectively. If $T_1$ is not a complete
binary tree, then $\psi(T)$ is defined to be the tree $T^*$ in
$\mathscr{P}_n\setminus (\mathscr{P}_n^*\cup \tilde{\mathscr{P}}_n)$
by replacing $T_1'$ in $T$ by $\psi(T_1')$; If $T_1$ is a complete
binary tree, so $T_2$ must not be a complete binary tree, then
$\psi(T)$ is defined to be the tree $T^*$ in $\mathscr{P}_n\setminus
(\mathscr{P}_n^*\cup \tilde{\mathscr{P}}_n)$ by replacing $T_2'$ in
$T$ by $\psi(T_2')$. See Figure \ref{fDD5} for example.
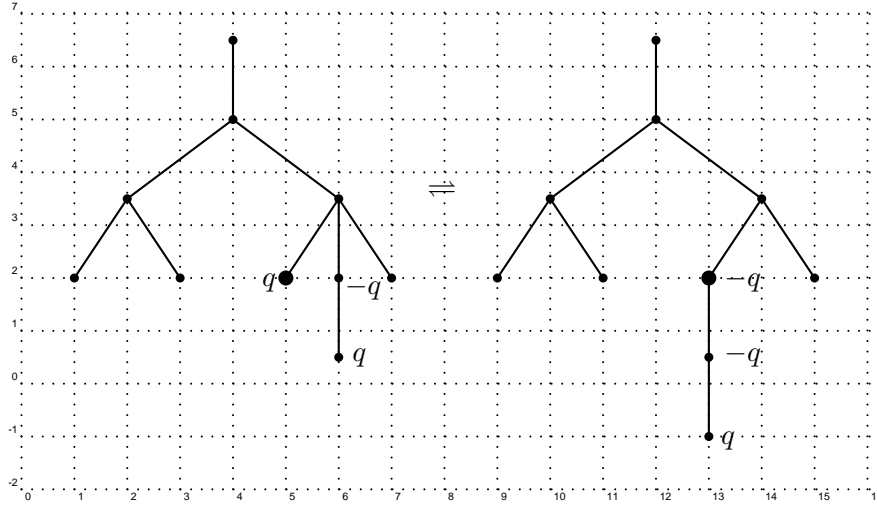
\begin{figure}[h]
\setlength{\unitlength}{0.5mm}
\begin{center}
\begin{pspicture}(11,4.8)
\psset{xunit=20pt,yunit=20pt}\psgrid[subgriddiv=1,griddots=5,
gridlabels=4pt](0,-2)(16,7)

\psline(4,6.5)(4,5)(2,3.5) \psline(4,5)(6,3.5)(5,2)
\psline(2,3.5)(1,2)\psline(2,3.5)(3,2)\psline(6,3.5)(6,2)(6,.5)\psline(6,3.5)(7,2)

\pscircle*(4,6.5){0.06}\pscircle*(4,5){0.06}
\pscircle*(2,3.5){0.06}\pscircle*(6,3.5){0.06}
\pscircle*(1,2){0.06}\pscircle*(3,2){0.06}\pscircle*(5,2){0.1}\pscircle*(6,2){0.06}
\pscircle*(6,.5){0.06}\pscircle*(7,2){0.06}

\put(4.3,1.2){$-q$}\put(4.4,.3){$q$}\put(3.2,1.3){$q$}
\put(5.4,2.5){$\rightleftharpoons$}

\psline(12,6.5)(12,5)(10,3.5)
\psline(12,5)(14,3.5)(13,2)(13,.5)(13,-1)
\psline(10,3.5)(9,2)\psline(10,3.5)(11,2)\psline(14,3.5)(15,2)

\pscircle*(12,6.5){0.06}\pscircle*(12,5){0.06}
\pscircle*(10,3.5){0.06}\pscircle*(14,3.5){0.06}
\pscircle*(9,2){0.06}\pscircle*(11,2){0.06}\pscircle*(13,2){0.1}\pscircle*(15,2){0.06}
\pscircle*(13,.5){0.06}\pscircle*(13,-1){0.06}

\put(9.35,1.3){$-q$}\put(9.3,-.8){$q$}\put(9.35,.3){$-q$}
\end{pspicture}\vskip1.3cm
\caption{The involution on $\mathscr{P}_n\setminus
\tilde{\mathscr{P}}_n$}\label{fDD5}
\end{center}
\end{figure}

\end{itemize}
Clearly, the $\psi$ as defined is indeed a sign reversing involution
on the set $\mathscr{P}_n\setminus \tilde{\mathscr{P}}_n$.
\end{proof}

Let $\mathscr{Q}_{n,k}$ denote the set of weighted plane trees of
$n+2$ vertices such that
\begin{itemize}
\item The leaves have weight $q^2$;
\item There exist $n-k$ vertices of out-degree one, except for the root, with
weight $-1$, $2q$ or $-q^2$, in other words, such vertices can be
regarded to be weighted by  $-(1-q)^2$;
\item All other internal vertices have weight $1$, there may exist
vertices of out-degree one with weight 1.
\end{itemize}
Let $\mathscr{Q}_n=\bigcup_{k=0}^n\mathscr{Q}_{n,k}$, and
$\mathscr{Q}_n^*$ denote the subset of $\mathscr{Q}_n$ such that
there are at least one vertex, except for the root, of out-degree
one weighted by $1$ or $-1$, let $\tilde{\mathscr{Q}}_n$ denote the
subset of $\mathscr{Q}_n$ such that the root has out-degree one, and
all other internal vertices have either out-degree two or out-degree
one with weight $2q$. It is clear that $\mathscr{Q}_n^*\cap
\tilde{\mathscr{Q}}_n=\emptyset$. Similar to the proof of Theorem
\ref{theo 4.1}, a sign reversing involution can be constructed on
$\mathscr{Q}_n\setminus \tilde{\mathscr{Q}}_n$, the detail leaves to
the interested readers. Hence we have
\begin{theorem}\label{theo 4.2}
There exists a sign reversing involution on the set
$\mathscr{Q}_n\setminus \tilde{\mathscr{Q}}_n$.
\end{theorem}

Now we can give the combinatorial proof of identities \eqref{eqn
3.8} and \eqref{eqn 3.9}.

\begin{proof}[Proof of \eqref{eqn
3.8} and \eqref{eqn 3.9}] For any $T\in \mathscr{P}_{k,k}$, namely,
$T$ is a weighted plane tree of $k+2$ vertices with leaves weighted
by $q$ and all other internal vertices weighted by $1$, inserting
$n-k$ vertices of weight $-(1+q)$ into the $k+1$ edges of $T$
(repetition allowed), we can obtain $\binom{n}{k}$ number of
weighted plane trees in $\mathscr{P}_{n,k}$. It is well-known that
the weight of $\mathscr{P}_{k,k}$ is $\mathfrak{N}_{k+1}(q)$, then
$\mathscr{P}_{n,k}$ has the weight
$\binom{n}{k}\mathfrak{N}_{k+1}(q)(-1-q)^{n-k}$. Similarly,
$\mathscr{Q}_{n,k}$ has the weight
$(-1)^{n-k}\binom{n}{k}\mathfrak{N}_{k+1}(q^2)(1-q)^{2(n-k)}$.

On the other hand, for any $T\in \tilde{\mathscr{P}}_n$, we know
that the root of $T$ has out-degree one and has only one subtree
$T'$ which is a weighted complete binary tree with $n+1$ vertices,
it is well known that the number of complete binary trees with $n+1$
vertices is counted by Catalan number $C_{\frac{n}{2}}$, where
$C_{\frac{n}{2}}=0$ if $n$ is odd. So the weight of
$\tilde{\mathscr{P}}_n$ is $q^{\frac{n}{2}+1}C_{\frac{n}{2}}$.

For any $T\in \tilde{\mathscr{Q}}_n$, let
$\tilde{\mathscr{Q}}_{n,k}$ denote the subset of
$\tilde{\mathscr{Q}}_n$ such that $T$ has $n-2k$ vertices, except
for the root, of out-degree one with weight $2q$. For any $T\in
\tilde{\mathscr{Q}}_{2k,k}$, we know that the root of $T$ has
out-degree one and has only one subtree $T'$ which is a weighted
complete binary tree with $2k+1$ vertices, inserting $n-2k$ vertices
of weight $2q$ into the $2k+1$ edges of $T$ (repetition allowed), we
can obtain $\binom{n}{2k}$ number of weighted plane trees in
$\tilde{\mathscr{Q}}_{n,k}$. It is clear that
$\tilde{\mathscr{Q}}_{2k,k}$ is counted by Catalan numbers $C_k$ and
has weight $q^{2k+2}C_k$, then $\tilde{\mathscr{Q}}_{n,k}$ has
weight
$(2q)^{n-2k}\binom{n}{2k}q^{2k+2}C_k=q^{n+2}\binom{n}{2k}C_k2^{n-2k}$.
Hence $\tilde{\mathscr{Q}}_{n}$ has weight
$q^{n+2}\sum_{k=0}^{n}\binom{n}{2k}C_k2^{n-2k}$, which is
$q^{n+2}C_{n+1}$ by Touchard identity.

Using Theorem \ref{theo 4.1} and \ref{theo 4.2}, one can easily
obtain that the weight of $\mathscr{P}_n ({\rm resp.}\
\mathscr{Q}_n)$ equals that of $\tilde{\mathscr{P}}_n ({\rm resp.}\
\tilde{\mathscr{Q}}_n)$, which completes the proof.
\end{proof} \vskip.5cm

\section*{Acknowledgements} The
second author was supported by The National Science Foundation of
China. 


\begin{thebibliography}{99}

\bibitem{aigner} M. Aigner, Catalan-like Numbers and Determinants,
{\em J. Combin. Theory, Series A} {\bf 87} (1999) 33--51.

\bibitem{bonin} J. Bonin, L. Shapiro and R. Simion, Some
$q$-analogues of the Schr\"{o}der numbers arising from combinatorial
statistics on lattice paths, {\em J. Statist. Plann. Inference} {\bf
34} (1993) 35--55.

\bibitem{chapman} R. Chapman, A curious identity, revisited, {\em The Math. Gazette}
{\bf 87} (2003) 139--141.

\bibitem{chendeng} W.Y.C. Chen, E.Y.P. Deng and R.R.X. Du,
Reduction of m-regular noncrossing partitions, {\em European J.
Combin.} {\bf 26}(2) (2005) 237--243.

\bibitem{chenpeng} W.Y.C. Chen and S. X. M. Pang, On the combinatorics of the Pfaff
identity, submitted.


\bibitem{chenyan} W.Y.C. Chen, S.H.F.Yan and L.L.M. Yang,
Identities from weighted $2$-Motzkin paths, {\em Adv. Appl. Math.},
to appear.

\bibitem{chenli} W.Y.C. Chen, N.Y. Li and L. Shapiro, The butterfly decomposition of plane
trees, {\em Disc. Appl. Math.}, {\bf 155}(17) 2187--2201.

\bibitem{coker} C. Coker, Enumerating a class of lattice paths, {\em Disc. Math.}
{\bf 271} (2003) 13--28.

\bibitem{comtet} L. Comtet, Advanced Combinatorics, D.Reidel, Dordrecht-Holland, 1970.

\bibitem{corsani} C. Corsani, D. Merlini and R. Sprugnoli, Left-inversion of
combinatorial sums, {\em Disc. Math.} {\bf 180}:1-3 (1998) 107--122.

\bibitem{deutsch}
E. Deutsch, Dyck path enumeration, {\em Disc. Math.} {\bf 204}
(1999) 167--202.

\bibitem{gasper} G. Gasper and M. Rahman, Basic Hypergeometric Series, 2nd Edition,
Encyclopedia of Mathematics And Its Applications, vol. 96, Cambridge
University Press, 2004.

\bibitem{hirschhorn} M. D. Hirschhorn, Comment on a ``curious"
identity, {\em The Math. Gazette} {\bf 87} (2003) 528--530.

\bibitem{mansun1} T. Mansour and Y. Sun, Dyck Paths and partial Bell polynomials,
{\em Australasian J. Combinatorics}., to appear.

\bibitem{mansun2} T. Mansour and Y. Sun, Bell polynomials and $k$-generalized Dyck Paths,
{\em Disc. Appl. Math.}, to appear.

\bibitem{munarini} E. Munarini, Generalization of a binomial
identity of Simons, {\em Integer: Elect. J. Combin. Number Theory}
{\bf 5} (2005) $\#$A15.

\bibitem{prodinger}
H. Prodinger, A curious identity proved by Cauchy's integral
formula, {\em The Math. Gazette} {\bf 89} (2005) 266--267.

\bibitem{riordan}
J. Riordan, Combinatorial Identities, New York: John Wiley \& Sons,
Inc. 1968

\bibitem{rogers} D. G. Rogers, Rhyming schemes: crossings and coverings, {\em
Disc. Math.} {\bf 33} (1981) 67--77.

\bibitem{rogershap} D. G. Rogers, L.W. Shapiro, Deques, trees and lattice paths,
{\em Lect. Notes Math.} {\bf 884} (1981) 293--303.

\bibitem{simull} R. Simion and D. Ullman, On the structure of the lattice of
noncrossing partitions, {\em Disc. Math.} {\bf 98} (1991) 193--206.

\bibitem{simons} S. Simons, A curious identity, {\em The Math. Gazette} {\bf 85} (2001) 296--298.

\bibitem{Slo}
N.J.A. Sloane, The On-Line Encyclopedia of Integer Sequences, {\em
www.research.att.com/jas/sequences}.

\bibitem{stanley} R. Stanley, Enumerative Combinatorics, vol. 1,
Cambridge Univ. Press, Cambridge, 1997.

\bibitem{sulanke} R. Sulanke, Counting lattice paths by Narayana polynomials, {\em
Electron. J. Combin.} {\bf 2} (2000) \#R40.

\bibitem{sun}
Y. Sun, The statistic ``number of $udu$'s" in Dyck paths, {\em Disc.
Math.} {\bf 287:1-3} (2004) 177-186.

\bibitem{wangsun}
X. Wang and Y. Sun, A new proof of a curious identity, {\em The
Math. Gazette} {\bf 91}:(502) (2007)  105--106.

\bibitem{wilf}  H. Wilf, Generatingfunctionology,
Academic Press, New York, 1990.
\end{thebibliography}
\end{document}